\documentclass{article}

\usepackage{amsmath, amsthm, amsfonts, amssymb, faktor, mathrsfs, textcomp, epsfig, graphicx, url}
\usepackage{mathrsfs}

\usepackage{ amsmath, amsfonts, amssymb, faktor, mathrsfs, graphicx}
\usepackage[curve, arrow]{xypic}

\newcommand{\op}{\operatorname}
\newcommand{\bull}{\mbox{\tiny{$\bullet$}}}

\theoremstyle{plain}
\newtheorem{thm}{Theorem}[section]
\newtheorem{prop}[thm]{Proposition}
\newtheorem{cor}[thm]{Corollary}
\newtheorem{lem}[thm]{Lemma}

\theoremstyle{definition}
\newtheorem{defin}[thm]{Definition}

\theoremstyle{remark}
\newtheorem{rem}[thm]{Remark}

\newtheorem*{notation}{Notation}

\begin{document}

\title{Cohomological Invariants of hyperelliptic curves of even genus}
\author{Roberto Pirisi}



\maketitle 

\begin{abstract}
Let $g$ be an even positive integer, and $p$ be a prime number. We compute the cohomological invariants with coefficients in $\mathbb{Z}/p\mathbb{Z}$ of the stacks of hyperelliptic curves $\mathscr{H}_g$ over an algebraically closed field $k_0$.
\end{abstract}

\maketitle 
\section{Introduction}
\begin{notation} We fix a base field $k_0$ and a prime number $p\neq \op{char}(k_0)$. All schemes and algebraic stacks will be assumed to be of finite type over $k_0$. If $X$ is a $k_0$-scheme we will denote by $\op{H}^{i}(X)$ the $i$-th \'etale cohomology group of $X$ with coefficients in $\mu_p^{\otimes i}$ (here $\mu_p^{\otimes 0}:=\mathbb{Z}/p\mathbb{Z}$), and by $\op{H}^{\bull}(X)$ the direct sum $\oplus_i \op{H}^i(X)$. If $R$ is a $k_0$-algebra, we set $\op{H}^{\bull}(R)=\op{H}^{\bull}(\op{Spec}(R))$.
\end{notation}

Cohomological invariants of algebraic groups have been widely studied in the last three decades, culminating in a number of important results by Merkurjev, Serre, Rost, and many others. They are invariants of the functor describing the isomorphism classes of $G$-torsors, which also naturally makes them invariants of the classifying stack $BG$.

 A natural extension of the concept is to construct corresponding invariants for more complicated moduli problems, such as for example smooth curves of a given genus. Differently from the case of $BG$, the classifying stacks for these moduli problems are not gerbes over a point. One could say that the stack $BG$ is a purely arithmetic object, which from a geometric point of view is just a point, with a group of automorphisms attached. On the opposite end of the spectrum, a variety can be thought as a purely geometric object. Stacks such as $\mathscr{M}_g$, which parametrizes families of smooth curves of genus $g$, have both geometric information, which we can roughly think about as the moduli space, and arithmetic information, corresponding to the curves' automorphism groups.

In \cite{Pirisi} the author introduced the notion of cohomological invariant of a smooth algebraic stack. Given a smooth algebraic stack $\mathscr{M}$, we can consider the functor of isomorphisms classes of its points 

\begin{center}
$P_{\mathscr{M}}: \left(\faktor{\textit{field}}{k_0}\right) \rightarrow \left( \textit{set} \right)$
\end{center}

which sends a field $K/k_0$ to the set of isomorphism classes of objects over $K$ in $\mathscr{M}$. Then a cohomological invariant for $\mathscr{M}$ is defined as a natural transformation

\begin{center}
$\alpha: P_{\mathscr{M}} \rightarrow \op{H}^{\bull}(-)$
\end{center}
satisfying a natural continuity condition.

The theory set up in \cite{Pirisi} can be used to compute the cohomological invariants of the stacks of hyperelliptic curves of all even genera. The main result of this paper is the following:

\begin{thm}
Suppose our base field $k_0$ is algebraically closed, of characteristic different from $2,3$. Let $g$ be an even positive integer, and let $\mathscr{H}_g$ be the stack of smooth hyperelliptic curves of genus $g$.
\begin{itemize}

\item For $p=2$, the cohomological invariants of $\mathscr{H}_{g}$ are freely generated as a graded $\mathbb{F}_2$-module by $1$ and invariants $\alpha_1,\ldots, \alpha_{g+2}$, where the degree of $\alpha_i$ is $i$.

\item If $p \neq 2$, then the cohomological invariants of $\mathscr{H}_{g}$ are nontrivial if and only if $2g + 1$ is divisible by $p$. In this case they are freely generated by $1$ and a single invariant of degree one.
 
\end{itemize}
\end{thm}

Moreover, we obtain partial results for non algebraically closed fields, proving in general that the cohomological invariants above are freely generated as a $\op{H}^{\bull}(k_0)$-module by the same elements if $p$ differs from $2$. If $p$ is equal to $2$, we show that the cohomological invariants of $\mathscr{M}_2$ can be decomposed in freely generated $\op{H}^{\bull}(k_0)$-module, whose generators are the same as in the algebraically closed case except for the one of highest degree, and a module $K$, which can be seen as a submodule of $\op{H}^{\bull}(k_0)$ shifted in degree by $4$.

The computation uses heavily Rost's theory of Chow groups with coefficients \cite{Rost}, and its equivariant version, which was first introduced by Guillot in \cite{Guillot}. This is due to the fact that for a smooth quotient stack $\left[ X /G \right]$ the zero-codimensional equivariant Chow group with coefficients $A^0_{G}(X,\op{H}^{\bull})$ is equal to the ring of cohomological invariants $\op{Inv}^{\bull}(\left[ X/G \right])$.

Our method is based on the presentation by Arsie and Vistoli \cite[4.7]{ArsieVistoli} of the stacks of hyperelliptic curves as the quotient of an open subset of an affine space by $\textit{GL}_2$. We use a technique similar to the stratification method introduced by Vezzosi in \cite{Vezzosi} and used by various authors afterwards (\cite{Guillot},\cite{MolinaVistoli}). The idea is, given a quotient stack $\left[ X/G \right]$, to construct an equivariant stratification $X=X_0\supset X_1 \ldots \supset X_n= \emptyset$ of $X$ such that the (equivariant) geometry of $X_i \! \smallsetminus \! X_{i+1}$ is simple enough that we can compute inductively the invariants for $X_i$ using the result for $X_{i+1}$ and the localization exact sequence \cite[p. 356]{Rost}.

One flaw of our method of computation is that it does not provide any real insight on the product in the ring of cohomological invariants of $\mathscr{H}_g$. The reason is that repeatedly applying the localization exact sequence causes loss of information about our elements, making it very difficult to understand what their products should be.

\section*{Description of content}

In section 2, we define a theory of Chern classes for Chow groups with coefficients, following \cite{ElmanKarpenkoMerkurjev}.  All the basic formulae for Chern classes hold with coefficients, and we use them to compute the Chow rings with coefficients of Grassmann bundles.

The second part of the section is dedicated to the equivariant version of the theory, which was introduced in \cite{Guillot}. We show that for a quotient stack the zero codimensional equivariant Chow group with coefficients computes the ring of cohomological invariants. To see that the whole theory translates to the equivariant setting we can just repeat the proofs used for the ordinary Chow groups in \cite{EdidinGraham}. We then compute the Chow rings with coefficients of $\textit{BGL}_n$ and $\textit{BSL}_n$.

In the third section we describe the presentation of the stacks $\mathscr{H}_g$, together with some related spaces that will appear in the computation. Consider $\mathbb{A}^{2g+3}$ as the space of binary forms of degree $2g+2$, and let $X_g \subset \mathbb{A}^{2g+3}$ be the open subscheme of nonzero forms with distinct roots. Then we have $\mathscr{H}_g=\left[ X_g / GL_{2} \right]$ for even $g$. The starting idea to compute the cohomological invariants of $\left[ X_g / GL_{2} \right]$ is that we can first compute those of $\left[ Z_g /\textit{GL}_2 \right]$, where $Z_g$ is the $\textit{GL}_2$-equivariant quotient $Z_g = X_g /G_m$. Here $G_m$ acts by multiplication. We introduce a stratification $P^{2g + 2} \supset \Delta_{1,2g+2} \supset \ldots \supset \Delta_{g+1,2g+2}$ which will be the base of our computation. We can see $\Delta_{i,2g+2}$ as the closed subscheme of binary forms divisible by the square of a form of degree $i$, so we have $Z_g = P^{2g+2} \! \smallsetminus \! \Delta_{1,2g+2}$.

In the fourth section we compute the invariants for $\mathscr{H}_g$ for all even $g$. The argument is based on the fact that the equivariant Chow groups with coefficients of $\Delta_{i,n} \! \smallsetminus \! \Delta_{i+1,n}$ are isomorphic to those of $P^{n-2i} \! \smallsetminus \! \Delta_{1,n-2i} \times P^{i}$, giving rise to an inductive reasoning relying heavily on the localization sequence. 

The fifth section is dedicated to extending the previous results to fields that are not algebraically closed. The extension turns out to be immediate when the prime $p$ is different from $2$, and rather troublesome for $p=2$. The main difficulty lies in understanding whether the pushforward of some elements is zero or not in the equivariant Chow ring with coefficients of $P^6$. To prove it, we construct an element in $A^{\bull}_{\textit{GL}_2}(P^6)$ that belongs to the annihilator of these images, but whose annihilator cannot contain them unless they are zero.

\emph{Aknowledgements: I wish to thank my Phd advisor Angelo Vistoli for making all of this possible. I wish to thank Burt Totaro and Zinovy Reichstein for their suggestions and careful reading of this material.}
\section{Chow rings with coefficients}

\subsection{Chern classes}
Chow groups with coefficients were first defined in Markus Rost's paper \cite{Rost}. They form a theory analogous to the classical theory of Chow groups, but with added arithmetic properties coming from the fact that the coefficients are taken, rather than in $\mathbb{Z}$, in a cycle module $M$. The latter is a functor from fields to graded abelian groups satisfying a long list of properties. The two most important examples of cycle modules are Milnor's $K$-theory and Galois Cohomology, in which case Chow groups with coefficients actually contain the ordinary Chow groups and the ordinary Chow groups mod $p$ respectively.

 Chow groups with coefficients have all the versatility of ordinary Chow groups, and moreover they satisfy a long exact sequence extending the short exact sequence of Chow groups. For a brief introduction to the subject the reader can refer to \cite[sec.2]{Guillot}. The entire theory is reworked in the greater generality of quasi-separated algebraic spaces in the author's PhD thesis \cite[ch.2]{Pirisi15}.

The \emph{Chow group with coefficients of dimension i}, written $A_{i}(X,M)$ is defined as the $i$-th homology group of the complex:

$$0 \rightarrow C_{d=\op{dim}(X)}(X,M) \rightarrow C_{d-1}(X,M) \rightarrow \ldots \rightarrow C_{1}(X,M) \rightarrow C_{0}(X,M) \rightarrow 0$$

Where $C_{i}(X,M)=\bigoplus_{P \in X_{(i)}} M(k(P))$ is the direct sum of the cycle module computed in the generic points of closed irreducible subschemes of $X$ of dimension $i$, and the differential comes from data attached to the cycle module \cite[pp.329,337]{Rost}. When each connected component of $X$ is equidimensional we can group points by codimension so that we get the complex:

$$0\rightarrow C^{0}(X,M) \rightarrow C^{1}(X,M) \rightarrow \ldots \rightarrow C^{\op{dim}(X)}(X,M) \rightarrow 0$$

With $C^{j}(X,M)=\bigoplus_{P\in X^{(j)}}M(k(P))$. Then the \emph{Chow group with coefficients of codimension i}, written $A^{i}(X,M)$, is defined as the $i$-th cohomology group of the complex above. When $X$ is moreover equidimensional we have the equality $C^j(X,M) =C_{\op{dim}(X)-j}(X,M)$ and consequently $A^{i}(X,M)=A_{\op{dim}(X)-i}(X,M)$.

If the scheme $X$ is smooth over $k_0$ and $M$ admits a bilinear pairing $M \times M \rightarrow M$ then we get a product $A^{i}(X,M)\times  A^{j}(X,M) \rightarrow A^{i+j}(X,M)$. In this case we call $A^{\bull}(X,M)$ the \emph{Chow ring with coefficients} of $X$. Note that it is a bigraded ring: we have the degree coming from the cycle module, and codimension.
 
Just as ordinary Chow groups, Chow groups with coefficients have a flat pullback, a proper pushforward, a pullback for maps to smooth schemes, they satisfy a projection formula and have a localization long exact sequence.
 
Rost's paper notably lacks the definition of a theory of Chern classes ``with coefficients". This has been done when $M$ is equal to Milnor's $K$-theory in chapter $9$ of Elman, Karpenko and Merkurjev's book \cite{ElmanKarpenkoMerkurjev}. We will extend their idea to all cycle modules. Our approach is slightly more cycle-based than the approach in \cite{ElmanKarpenkoMerkurjev}.

\begin{defin}
Let $E\rightarrow X$ be a vector bundle. A \emph{coordination} $\sigma$ for $E$ is a sequence of closed subsets $(X=X_0 \supset X_1 \supset \ldots \supset X_n)$ together with trivializations of $E$ along each of the locally closed subschemes $X_0 \smallsetminus X_1, \ldots ,X_{n-1}\smallsetminus X_{n}, X_{n} $. 
\end{defin}

\begin{defin}\label{C1def}
Let $L \xrightarrow{\pi} X$ be a line bundle. Let $\sigma$ be a coordination for $L$ and let $i$ be the zero-section imbedding. In \cite[sec. 9]{Rost} Rost defines a retraction $r_{\sigma}:C^{\bull}(L,M) \rightarrow C^{\bull}(X,M)$, depending on the coordination $\sigma$, for the pullback $C^{\bull}(X,M) \rightarrow C^{\bull}(L,M)$.  We define the \emph{first Chern class} $c_{1,\sigma}(L): C_p(X,M) \rightarrow C_{p-1}(X,M)$ of the couple $(L,\sigma)$ as 
$$c_{1,\sigma}(L)(\alpha) = r_{\sigma} \circ i_* (\alpha)$$
\end{defin}

Clearly the choice of a coordination is irrelevant in homology and we will just refer to $c_1(L)$ when we are interested in the induced map in homology. The additional data of the coordinations allows for slightly more precise statements on cycle level when we pull back the coordination together with the line bundle, as we will see:

\begin{prop}\label{C1}
Consider a morphism $Y \xrightarrow{f} X$ and form the cartesian square:
\begin{center}
$\xymatrixcolsep{5pc}
\xymatrix{ E \ar@{->}[r]^{f_1} \ar@{->}[d]^{\pi_1} & L \ar@{->}[d]^{\pi} \\ 
Y \ar@{->}[r]^{f} & X  }$
\end{center}
Let $\sigma'$ be the induced coordination on $Y$. Then:
\begin{enumerate}
\item If $f$ is proper then $f_*( c_{1,\sigma'}(L)(\alpha)) = c_{1,\sigma}(L)( f_*(\alpha))$.
\item If $f$ is flat then $c_{1,\sigma'}(L)(f^*(\beta))=f^*(c_{1,\sigma}(L)(\beta))$.
\item If $a \in \mathcal{O}_X^*(X)$ then $c_{1,\sigma}(\lbrace a \rbrace (\alpha))=\lbrace a \rbrace (c_{1,\sigma}(\alpha))$.
\item Suppose $X$ is normal. Suppose $E=\mathcal{O}(D)$ for an irreducible subvariety $D$ of codimension $1$, defined by a valuation $v$ on $k(X)$, and let $\sigma$ be a coordination with $X_1=D$. Then the restriction of $c_{1,\sigma}(E)$ to $A^0(X)$ is equal to the map $s_v^{\tau}$ defined in \cite[R3f]{Rost}. In particular, $c_1(E)(1) = D$.

The following properties are true at homotopy level:

\item If $(U,V)$ is a boundary couple then $\partial^U_V(c_1(L_{\mid U})(\alpha))=c_1(L_{\mid V})(\partial^U_V(\alpha))$.
\item If $L,E$ are two line bundles over $X$ then $c_1(L)(c_1(E)(\alpha))=c_1(E)(c_1(L)(\alpha))$.
\item The first Chern class of $L\otimes E$ satisfies $c_1(L\otimes E)(\alpha) = c_1(L)(\alpha) + c_1(E)(\alpha)$.
\end{enumerate}
\end{prop}
\begin{proof}
Properties (i),(ii),(iii) can be immediately obtained by the compatibility of the retraction $r_{\sigma}$ with the fundamental maps \cite[9.5]{Rost}. 

 We prove property (iv) by an explicit local computation. Given an element $\alpha \in A^0(X)$ we can explicitly write down $c_{1,\sigma}(L)(\alpha)$. It is equal to $r_{\sigma_{\mid D}} \circ \partial^{L_{X\! \smallsetminus \! D}}_{L_D} \circ H_{\textit{triv}} \circ (i_0)_* (\alpha)$. 

By explicit verification we see that $H_{\textit{triv}} \circ (i_0)_* (\alpha)=\lbrace t \rbrace( \pi^* (\alpha))$, where $t$ is the parameter for the trivial bundle over $X \! \smallsetminus \! D$. The expression makes sense as the cycle $\pi^* (\alpha)$ is not supported in any point where $t$ is zero.

Now we consider the boundary map $\partial^{L_{X\! \smallsetminus \! D}}_{L_D}$. As $\mu = \lbrace t \rbrace( \pi^* (\alpha))$ lives in the generic point of $L$, the only point where the value of $\partial^{L_{X\! \smallsetminus \! D}}_{L_D}(\mu)$ can be nonzero is the generic point of $L_D$. Then $\partial^{L_{X\! \smallsetminus \! D}}_{L_D}(\mu)$ is equal to the value of $\mu$ through the map $\partial_{v}:\op{H}^{\bull}(k_0(L))\rightarrow \op{H}^{\bull}(k_0(L_D))$. To compute it we first base change $\mu$ to a neighborhood $U_D \xrightarrow{\rho} X$ of the generic point of $D$ such that the bundle is trivial.

 In base changing $\mu$ to $U_D$ we need to keep track of what happens to $\lbrace t \rbrace$, which is no longer the parameter for our trivial bundle: if $t'$ is the new parameter, we see that we can write $t= \tau t'$, where $\tau \in \mathcal{O}^*_X(U_D \times_X (X \! \smallsetminus \! D))$ vanishes in $D$ with order $1$. Then the pullback of $\mu$ to $U_D$ is equal to $\lbrace t' \rbrace  (\pi^* \alpha) + \lbrace \tau \rbrace  (\pi^* \alpha)= \mu_1 + \mu_2 $. 
 
 Now it's easy to see that $$\partial^{L_{X\! \smallsetminus \! D}}_{L_D}(\mu_1)=\lbrace t' \rbrace\pi^*\partial_D(\alpha)=0$$ and $$\partial^{L_{X\! \smallsetminus \! D}}_{L_D}(\mu_2)=\partial_{v}(\tau \pi^* \alpha)=s^{\tau}_{v_L}(\pi^*(\alpha))$$ where $v_L$ is the valuation defining $L_D$. By the compatibility of the map $s^{\tau}_{v_L}$ with the retraction $r_{\sigma}$ we obtain the result. It should be noted that at the homology level the maps $s^{\tau}_v$ do not depend on the choice of the uniformizing parameter $\tau$.
 
 Property (v) is obtained using the compatibility of pullback and differential, by writing $r_{\sigma}=\pi^{-1}$ as homology maps.
 
We still need to prove the last two properties.

Consider a cartesian square:
\begin{center}
$\xymatrixcolsep{5pc}
\xymatrix{ L\times_X E \ar@{->}[r]^{\pi'_1} \ar@{->}[d]^{\pi_1} & E \ar@{->}[d]^{\pi} \\ 
L \ar@{->}[r]^{\pi'} & X  }$
\end{center}
And name $i,i'$ the zero sections respectively of $E$ and $L$, and $i_1,{i'}_1$ the zero sections of respectively $\pi$ and $\pi'$.

 Then $c_1(L)(c_1(E)(\alpha))=({\pi'^*})^{-1} \circ {i'}_* \circ (\pi^*)^{-1} \circ i_*$. By the compatibility with proper pushforward we have $({\pi'^*})^{-1} \circ {i'}_* \circ (\pi^*)^{-1} \circ i_* = ({\pi'^*})^{-1} \circ ({\pi'}^*_1)^{-1} \circ i_* \circ ({i'_1})_*$.
 
  By the functoriality of pullback and pushforward we get the equality $({\pi'^*})^{-1} \circ ({\pi'}^*_1)^{-1} \circ i_* \circ ({i'_1})_*= ((\pi' \circ \pi_1)^*)^{-1} \circ (i'_1 \circ i)_*$ which is equal to $((\pi'_1 \circ \pi)^*)^{-1} \circ (i_1 \circ i')_*$ as the two maps are the same, and doing the reasoning above backwards we obtain the desired equality (vi).
  
  For the last equality, recall that there is a flat product map $L\times_X E \xrightarrow{\rho} L \otimes E$ such that the composition of $\rho$ and the projection $\pi'':L\otimes E \rightarrow X$ is the projection $\pi_2: L\times_X E \rightarrow X$.
  
     It is easy to see that if $i''$ is the zero section of $\pi''$ then $\rho^* \circ i''_* (\alpha) = \pi_1^* \circ i'_* (\alpha) + {\pi'}_1^* \circ i_* (\alpha)$. As the projections from $E\times_X L$ and $E\otimes L$ to $X$ both induce an isomorphism we know that $\rho^*$ must be an isomorphism too. But then 
     
     $$c_1(E \otimes L)(\alpha) = ({\pi''}^*)^{-1} \circ i''_* (\alpha)= ({\pi''}^*)^{-1} \circ ((\rho^*)^{-1} \circ \rho^*)  \circ i'' =(\pi_2^*)^{-1} \circ \rho^* \circ i''$$
     Which is in turn equal to 
     $$ (\pi_2^*)^{-1}(\pi_1^* \circ i'_* (\alpha) + {\pi'}_1^* \circ i_* (\alpha))=c_1(L)(\alpha) + c_1(E)(\alpha)$$.
\end{proof}

As the maps we defined commute on homology level, we can see the composition $c_1(L_1) \circ c_1(L_2)$ as a commutative product $c_1(L_1) \circ c_1(L_2)= c_1(L_1)\cdot c_1(L_2)$.

Having this definition down we can follow closely the sections $53-55$ of \cite{ElmanKarpenkoMerkurjev} to obtain a complete theory of Chern classes. This is done in detail in the author's PhD thesis \cite[chapter 2]{Pirisi15}.

\begin{prop}\label{proj}
Let $P(E) \xrightarrow{\pi} X$ be the projectivization of a vector bundle $E$ of rank $r+1$ over $X$. The following formula holds for all $p$:
$$ A_{p}(E) = \bigoplus_{\substack{n-i=p \\ i \leq r}}  c_1(\mathcal{O}_{P(E)}(-1))^i (\pi^*(A_n (X,M))) \simeq \bigoplus_{\substack{p+r\leq n \leq p}} A_n (X,M) $$
\end{prop}
\begin{proof}
The proof in \cite[sec.53]{ElmanKarpenkoMerkurjev} holds without changes.
\end{proof}

\begin{prop}\label{Chern}
There is a theory of Chern classes for Chow groups with coefficients, satisfying all of the usual properties.
\end{prop}
\begin{proof}
This is done in the author's Phd thesis \cite[ch.2, sec.4]{Pirisi15}, and it can also be obtained easily following sections $53-55$ in \cite{ElmanKarpenkoMerkurjev}.
\end{proof}

For the rest of the section we assume that our cycle module $M$ has a pairing $M\times M\rightarrow M$, so that for smooth schemes we have a ring structure. With the next corollary we we show that for smooth schemes Chern classes are represented by multiplication by elements in the Chow ring with coefficients.

\begin{cor}\label{mult}
If $X$ is smooth over $k_0$, and $D$ is an irreducible divisor then $c_1(\mathcal{O}(D))(\alpha)=\alpha \cdot D$ and the classes $c_i(E)$ are all equal to the product by $\beta$ for some cycle $\beta$ of degree zero.
\end{cor}
\begin{proof}
Let $L=\mathcal{O}(D)$. As $L$ is smooth we can consider the pullback through the zero section $i: X \rightarrow L$ and by the compatibility with flat pullback we must have $i^*\circ\pi^*=\op{Id}_{A^{\bull}(X)}$. Now consider a product $\alpha \cdot \beta$. We can see $\alpha$ as $i^* \circ \pi^* \alpha$, so by the projection formula we have $i_*(\alpha \cdot \beta)=\pi^* \alpha \cdot i_* \beta$. Then as $c_1(L)=(\pi^*)^{-1}\circ i_*$ we have $c_1(L)(\alpha \cdot \beta)=(\pi^*)^{-1}(\pi^* \alpha \cdot i_* \beta)=\alpha \cdot c_1(L)(\beta)$. Then we can take $\beta = 1$ to get the identity $c_1(L)(\alpha)=\alpha \cdot c_1(L)(1)$ and as $c_1(L)(1)=D$ by \ref{C1} we can conclude.

 The general case is a direct consequence of the line bundle case by using the splitting principle and the Whitney sum formula.
\end{proof}

With the next proposition we use the theory available to us to compute the Chow groups with coefficients of Grassmann bundles. This will be instrumental in computing the equivariant Chow groups with coefficients of $\textit{GL}_2$ in the next section.
%

\begin{cor}\label{Grass}
Suppose $X$ is smooth over $k_0$, and $M$ is equal to either Milnor's $K$-theory or $M=\op{H}^{\bull}$. Let $E \xrightarrow{\pi} X$ be a vector bundle, and let $\op{Gr}_m(E) \rightarrow X$ be the Grassmann bundle of $m$-dimensional subbundles of $E$. Then $$A^{\bull}(\op{Gr}_m(E),M)= A^{\bull}(X,M)\otimes_{\op{CH}^{\bull}(X)}\op{CH}^{\bull}(\op{Gr}_m(E))$$
\end{cor}
\begin{proof}
%

We have an obvious map $$A^{\bull}(X,M)\otimes_{\op{CH}^{\bull}(X)}\op{CH}^{\bull}(\op{Gr}_m(E))\xrightarrow{\psi} A^{\bull}(\op{Gr}_m(E),M)$$ given by pulling back the first factor to $A^{\bull}(\op{Gr}_m(E))$ and multiplying
$$\alpha \otimes \lambda \xrightarrow{\psi} \alpha \cdot \lambda $$

Note now that we can obtain the bundle of complete flags $\op{FL}_m(E)$ from either $X$ or $\op{Gr}_m(E)$ by a sequence of projective bundles, via the usual splitting construction. Using propositions \ref{proj}, \ref{Chern} we can conclude that composing with the pullback to $\op{FL}_m(E)$ we obtain an injective map $$A^{\bull}(X,M)\otimes_{\op{CH}^{\bull}(X)}\op{CH}^{\bull}(\op{Gr}_m(E))\rightarrow A^{\bull}(\op{Fl}_m(E),M)$$ 

This shows injectivity. Let $\pi:A^{\bull}(\op{Fl}_m(E)) \rightarrow A^{\bull}(\op{Gr}_m(E))$ be the projection. We can construct a section to $\pi^*$ by composing the sections for each projective bundle: it amounts to multiplying by a fixed element of degree zero $x$ and then taking the pushforward $\pi_*$. Now consider an element $\beta \in A^{\bull}(\op{Gr}_m(E))$. Inside $A^{\bull}(\op{Fl}_m(E))$ we have $\alpha=\beta \cdot y$, where $\beta$ comes from $A^{\bull}(X)$ and $y$ has degree zero. Then $\alpha= \pi_* (\pi^*\alpha \cdot x)=\pi_*(\pi^*\beta \cdot x \cdot y)=\beta \cdot \pi_*(x \cdot y)$, which belongs to the image of $\psi$.
\end{proof}


\subsection{Equivariant theory}

A cycle-based approach is clearly only feasible when points have a well defined underlying field. Defining a theory of Chow groups with coefficients for a suitably large class of algebraic stacks will require a different approach. For a quotient stack $\left[ X/G \right]$ we can use the same type of equivariant approach defined in \cite{Totaro} and \cite{EdidinGraham}. This has already been described in \cite{Guillot}.

The basic idea is that any extension of the theory should be homotopy invariant, and that the Chow groups with coefficients of dimension $i$ should not change if we remove or modify somehow a subset of dimension at least $i-2$. Using this, up to readjusting the dimension index, we can replace our object of study $\left[ X/G \right]$ with a scheme $E \rightarrow \left[ X/G \right]$ that, up to some high codimension subset, is a vector bundle over $\left[ X/G \right]$.

In the next definition we use the fact that in the author's Phd thesis \cite{Pirisi15} a theory of Chow groups with coefficients for algebraic spaces is defined. One can add some technical conditions (e.g. quasi-projective scheme with a linearized group action, see \cite[prop.23]{EdidinGraham}), which will be true in every case we consider in this paper, to make sure all quotients involved are schemes.

\begin{defin}\label{equiv1}
Let $i$ be a positive integer, and let $X$ be a scheme equipped with an action by an algebraic group $G$. Let $V$ be a $r$-dimensional representation of $G$ such that $G$ acts freely outside of a closed subset $W=V\! \smallsetminus \! U$ of codimension equal or greater than $i+2$, and set $U=V\! \smallsetminus \! W$.

Consider the quotient $X \times^G U = (X \times U)/G$, where the action of $G$ is the diagonal action. By \cite[02Z2]{StacksProject} we know that $X \times^G U$ is an algebraic space, and if $X$ is quasi separated so is $X \times^G U$. In this case we define the \emph{equivariant Chow group with coefficients of dimension i} $A_i^G(X,M)$ to be $A_{i+r-\op{dim}(G)}(X \times^G U,M)$.

If $X$ is equidimensional can also switch to the codimensional notation by writing 
$$A^i_G(X)=A^i(X\times^G U,M)=A_{r+\op{dim}(X)-\op{dim}(G)-i}^G(X,M)$$
\end{defin} 

The groups above are well defined (i.e., they do not depend on the particular representation $V$) by the double fibration argument, as in \cite[2.2]{EdidinGraham}. The argument consists of the fact that if we have two representations $V,V'$ of dimension $r,r'$ satisfying the requirements for the definition we can construct a third representation $V \times V'$ and then $A_{i+r+r'}(X\times^G (U\times U'))$ is isomorphic to both $A_{i+r}(X\times^G U)$ and $A_{i+r'}(X\times^G U')$.

Note that there is no reason why the equivariant groups should be zero for codimension $\gg 0$ or negative dimension, and in fact this is not the case even for the most basic examples.

\begin{prop}
The equivariant Chow groups with coefficients $A^i_G(X,M)$ only depend on the quotient stack $\left[ X/G \right]$.

 Every result in \cite{Rost} can be restated for $G$-equivariant Chow groups with coefficients and $G$-equivariant maps.
\end{prop}
\begin{proof}
We can use the double fibration argument and the fact that equivariant maps are well-behaved when passing to equivariant approximation, as in \cite[5.2, prop.16]{EdidinGraham}, \cite[2.2, prop.3]{EdidinGraham}.
\end{proof}

\begin{prop}
Let $\left[X/G\right]$ be a quotient stack, smooth over $k_0$. Then $$A^0_G(X,\op{H}^{\bull})=\op{Inv}^{\bull}(\left[X/G\right]).$$
\end{prop}
\begin{proof}
Consider an equivariant approximation $\left[ X \times U /G \right] \xrightarrow{\pi} \left[ X/G \right]$ such that $T := \left[ X \times U /G \right]$ is an algebraic space. The pullback $\pi^*$ induces an isomorphism on cohomological invariants. This is an immediate consequence of the fact that cohomological invariants do not change when passing to a vector bundle or when removing a subset of codimension two or more \cite[4.13,4.14]{Pirisi}.

Now recall that on a scheme, or more generally on an algebraic space $T$ we know that $A^0(T,\op{H}^{\bull})=\op{Inv}^{\bull}(T)$ by \cite[4.9]{Pirisi}. Then we have a natural identification
$$ \op{Inv}^{\bull}(\left[X/G\right]) = \op{Inv}^{\bull}(T) = A^0(T,\op{H}^{\bull}) = A^0_G(X,\op{H}^{\bull})$$
\end{proof}

We will now compute some equivariant Chow groups with coefficients, taking $G$ a classical group acting trivially on the spectrum of a field. The computations for $G=\textit{GL}_n, \textit{SL}_n$ are an immediate consequence of our description of the Chow ring with coefficients for Grassmanian bundles.

Recall that the $\textit{GL}_n$-equivariant Chow ring of a point is

$$ \op{CH}^{\bull}_{\textit{GL}_n}(\op{Spec}(k_0))=\mathbb{Z}\left[ c_1,\ldots, c_n \right]$$

Where $c_1,\ldots,c_n$ are the Chern classes of the canonical representation. The $\textit{SL}_n$-equivariant Chow ring of a point is

$$ \op{CH}^{\bull}_{\textit{SL}_n}(\op{Spec}(k_0))=\mathbb{Z}\left[ c_2,\ldots, c_n \right]$$

Where again the $c_i$ are the Chern classes of the canonical representation.

\begin{prop}\label{GLSL}
Let our cycle module $M$ be either equal to Milnor's $K$-theory or Galois cohomology, and let $G$ be the general linear group $\textit{GL}_n$ or the special linear group $\textit{SL}_n$. Then the equivariant Chow ring with coefficients $A^{\bull}_G(\op{Spec}(k_0),M)$ is equal to the tensor product of the corresponding ordinary equivariant Chow groups and $M(\op{Spec}(k_0))$.
\end{prop}
\begin{proof}
As in \cite[3.1]{EdidinGraham}, by choosing an appropriate representation of $\textit{GL}_n$ we can compute its equivariant Chow groups with coefficients on Grassmann schemes. Then the description given in \ref{Grass} of the Chow groups of Grassmann bundles allows us to conclude immediately. Note that the product structure is the same as that induced on the tensor product by the two products on ordinary Chow groups and on $M(\op{Spec}(k_0))$; this can be seen as a consequence of proposition \ref{mult}.

 Now, following \cite[sec.3]{MolinaVistoli} if we consider the representation of $\textit{SL}_n$ induced by our original representation of $\textit{GL}_n$, we see that the natural map $\op{Spec}(k_0) \times^{\textit{SL}_n} U \rightarrow \op{Spec}(k_0) \times^{\textit{GL}_n} U$ is a $G_m$-torsor with associated line bundle the determinant bundle. Using the injectivity of $c_1$, a simple long exact sequence argument shows that the ring $A^{\bull}_{\textit{SL}_n}(\op{Spec}(k_0))$ must be isomorphic to $\faktor{A^{\bull}_{\textit{GL}_n}(\op{Spec}(k_0),M)}{c_1}$
\end{proof}

 \section{Preliminaries}

In this section we state some general considerations that will be needed for all the computations in the paper. From now on whenever we use Chow groups with coefficients (resp. equivariant Chow groups with coefficients) our cycle module will be $\op{H}^{\bull}$, so we shorten $A^{\bull}(X,\op{H}^{\bull})$ to $A^{\bull}(X)$ (resp. $A_{G}^{\bull}(X,\op{H}^{\bull})$ to $A_{G}^{\bull}(X)$).

We begin by recalling the presentations of the stacks we will work with, all due to Arsie and Vistoli \cite{ArsieVistoli}.

\begin{thm}\label{Pres}
Let $g$ be an even positive integer. Consider the affine space $\mathbb{A}^{2g + 3}$, seen as the space of all binary forms $\phi(S)=\phi(s_0,s_1)$ of degree $2g + 2$. Denote by $X$ the open subset consisting of nonzero forms with distinct roots, with the action of $\textit{GL}_2$ on $X_g$ defined by $A(\phi(S))=\op{det}(A)^{g} \phi(A^{-1}S)$.

 Denote by $\mathscr{H}_g$ the stack of smooth hyperelliptic curves of genus $g$. In particular, as any smooth curve of genus $2$ is hyperelliptic, $\mathscr{H}_2=\mathscr{M}_2$. Then we have 
 $$\mathscr{H}_{g} \simeq \left[ X_g/\textit{GL}_2\right]$$
And the canonical representation of $\textit{GL}_2$ yields the Hodge bundle of $\mathscr{H}_{g}$.

\end{thm}
\begin{proof}
This is corollary 4.7 of \cite{ArsieVistoli}. When $g=2$, the presentation of $\mathscr{M}_2$ was originally shown by Vistoli in \cite[3.1]{Vistoli}.
\end{proof}

In both cases, the quotient of $X_g$ by the usual action of $\mathbb{G}_m$ defined by $(x_0,\ldots,x_{2g+2},t) \rightarrow (tx_{0},\ldots,tx_{2g+2})$, which we will name $Z_g$, is naturally an open subset of the $\textit{GL}_2$-scheme  $P(\mathbb{A}^{2g + 3})$, namely the complement of the discriminant locus.

We will first construct the invariants of the quotient stack $\left[ Z_g/\textit{GL}_2 \right]$, then use the principal $\mathbb{G}_m$-bundle $\left[ X_g/\textit{GL}_2 \right] \rightarrow \left[ Z_g/\textit{GL}_2 \right]$ to compute the invariants of $\mathscr{H}_g$ for $g$ even.

We generalize the family of equivariant schemes in above this way: let $F$ be the dual of the standard representation of $\textit{GL}_2$. We can see $F$ as the space of all binary forms $\phi=\phi(s_0,s_1)$ of degree $1$. It has the natural action of $\textit{GL}_2$ defined by $A(\phi)(x)=\phi(A^{-1}(S))$. We denote by $E_i$ the $i$-th symmetric power $\op{Sym}^i(F)$. We can see $E_i$ as the space of all binary forms of degree $i$, with an action of $\textit{GL}_2$ induced by the action on $F$ which is $A(\phi)(S)=\phi(A^{-1}(S))$.

 Note that while this action is different from the action of $\textit{GL}_2$ on $X_g\subset E_{2g+2}$, it is the same when we pass to the projectivized schemes $Z_g$.

\begin{defin} We denote $\tilde{\Delta}_{r,i}$ the closed subspace of $E_i$ composed of forms $\phi$ such that there exists a form $f$ of degree $r$ whose square divides $\phi$. We can think of it as the image of the map
$$E_{i-2r} \times E_r \rightarrow E_i, \quad (f,g) \rightarrow fg^2 $$
which is a closed subset as the map passes to the projectivizations. With this notation the scheme $X_g$ in theorem \ref{Pres} is equal to $E_{2g+2} \! \smallsetminus \! \tilde{\Delta}_{1,2g+2}$. 

 We denote $\Delta_{r,i}$ the closed locus in the projectivizations $P(E_i)$ composed of forms $\phi$ such that there exists a form $f$ of degree $r$ whose square divides $\phi$. With this notation we have $Z_g = P(E_{2g+2}) \! \smallsetminus \! \Delta_{1,2g+2}$.
\end{defin}
Thanks to the localization exact sequence, understanding the cohomological invariants of $\left[ P(E_i) \! \smallsetminus \! \Delta_{1,i}/\textit{GL}_2 \right]$ can be reduced to understanding the invariants of $\left[ P(E_i) /\textit{GL}_2 \right]$, which are understood thanks to the projective bundle formula, the top Chow group with coefficients $A^0_{\textit{GL}_2}(\Delta_{1,i})$ (which is not equal to the group of cohomological invariants of $\left[ \Delta_{1,r}/\textit{GL}_2 \right]$, as $\Delta_{1,i}$ is not smooth) and the pushforward map $A^0_G(\Delta_{1,i})\rightarrow A^1_{\textit{GL}_2}(P(E_i))$. The computation of $A^0_{\textit{GL}_2}(\Delta_{1,i})$ will be based on the following two propositions.

\begin{prop}\label{univ}
Let $\pi_{r,i}: P(E_{i-2r})\times P(E_{r})\rightarrow \Delta_{r,i}$ be the map induced by $(f,g)\rightarrow fg^2$. The equivariant morphism $\pi_{r,i}$ restricts to a universal homeomorphism on $\Delta_{r,i} \! \smallsetminus \! \Delta_{r+1,i}$. Moreover, if $\op{char}(k_0) > 2r$ or $\op{char}(k)=0$ then any $k$-valued point of $\Delta_{r,i} \! \smallsetminus \! \Delta_{r,i+1}$ can be lifted to a $k$-valued point of $P(E_{i-2r})\times P(E_{r})$.
\end{prop}
\begin{proof}
See \cite[3.2]{Vistoli}. The reasoning holds in general as long as we can say that a polynomial with $r$ double roots must be divisible by the square of a polynomial of degree $r$. This is clearly true for $\op{char}(k)=0$, but in positive characteristic it holds only as long as $2r < \op{char}(k)$, as we can find irreducible polynomials of degree $\op{char}(k)$ with only one distinct root. It is however always true that the map $\pi_{r,i}$ is a bijection when restricted to $\Delta_{r,i} \! \smallsetminus \! \Delta_{r+1,i}$. Being proper and bijective, it is a universal homeomorphism.
\end{proof}

\begin{prop}\label{univiso}
The pushforward of a (equivariant) universal homeomorphism induces an isomorphism on (equivariant) Chow groups with coefficients in $\op{H}^{\bull}$.
\end{prop}
\begin{proof}
Note first that the non-equivariant statement implies the equivariant one, as if $X,Y$ are $G$-schemes on which $G$ acts freely then an equivariant universal homeomorphism between them induces a universal homeomorphism on quotients.

Let $f: X \rightarrow Y$ be a universal homeomorphism. Given a point $y \in Y$, its fibre $x$ is a point of $X$ and the map $f_x: x \rightarrow y $ is a purely inseparable field extension. The pullback $(f_x)^{*}: \op{H}^{\bull} (y) \rightarrow \op{H}^{\bull} (x)$ is an isomorphism, and the projection formula yields $(f_x )_* ( (f_x )^* \alpha ) = \left[ k(x) : k(y) \right] \alpha$. As the characteristic of $k(x)$ is different from $p$, the degree $\left[ k(x) : k(y) \right]$ is invertible modulo $p$ and the corestriction map is an isomorphism. This implies that $f_*$ induces an isomorphism on cycle level.

\end{proof}

In the following we will write $P^i$ for $P(E_i)$, which should not cause any confusion as $\op{dim}(E_i)=i+1$. We will exploit the stratification $$\Delta_{1,i} = \Delta_{1,i} \! \smallsetminus \! \Delta_{2,i} \sqcup \Delta_{2,i} \! \smallsetminus \! \Delta_{3,i} \sqcup \ldots \sqcup \Delta_{\left[ i/2 \right],i}$$ and the isomorphism $$A^0_G(\Delta_{r,i} \! \smallsetminus \! \Delta_{r+1,i}) \simeq A^0_{\textit{GL}_2}((P^{i-2r}\! \smallsetminus \! \Delta_{1,i-2r})\times P^r)$$ to inductively compute $A^0_{\textit{GL}_2}(\Delta_{1,i})$ and $A^0_{\textit{GL}_2}(P^i \! \smallsetminus \! \Delta_{1,i})$.

\section{The invariants of $\mathscr{H}_{g}$, $g$ even}

We are now going to compute the cohomological invariant of the stacks $\mathscr{H}_g$ for all even $g$. We assume that our base field $k_0$ is algebraically closed.

 Recall that as we are working with schemes over $k_0$ all the cohomological invariants and Chow groups with coefficients that we are going to consider are $\op{H}^{\bull}(k_0)$-modules, and all the maps are maps of $\op{H}^{\bull}(k_0)$ modules. Due to the fact that $k_0$ is algebraically closed, we have $\op{H}^{\bull}(k_0)=\op{H}^{0}(k_0)=\mathbb{Z}/p\mathbb{Z}$. Our main result is the following:

\begin{thm}\label{evenmain}
Suppose our base field $k_0$ is algebraically closed, of characteristic different from $2,3$, and let $g$ be even. For $p=2$, the cohomological invariants of $\mathscr{H}_{g}$ are freely generated as a graded $\mathbb{F}_2$-module by $1$ and invariants $\alpha_1,\ldots, \alpha_{g+2}$, where the degree of $\alpha_i$ is $i$. If $p \neq 2$, then the cohomological invariants of $\mathscr{H}_{g}$ are nontrivial if and only if $2g + 1$ is divisible by $p$. In this case they are generated by $1$ and a single nonzero invariant in degree one.
\end{thm}

A few last considerations on equivariant Chow rings are needed.

\begin{lem}\label{Diag}
Let $F$ be a vector bundle of rank $2$ on a variety $S$ smooth over $k_0$, let $P = P(F)$ be the projective bundle of lines in $F$, and $\Delta$ the image of the diagonal embedding $\delta : P \rightarrow P \times_{S} P$. Let $y_1, y_2$ in $\op{CH}^{\bull}(P\times_{S}P)$ be the two pullbacks of the first Chern class of $\mathcal{O}_{P}(1)$, $c_1 \in \op{CH}^{\bull}(P\times_{S}P)$ the pullback of the first Chern class of $F$. Then the class of $\Delta$ is $y_1 + y_2 +c_1$.
\end{lem}
\begin{proof}
This is \cite[3.8]{Vistoli}.
\end{proof}

Using the previous lemma we are able to compute the classes of $\Delta_{1,i}$ in $\op{CH}^{1}_{\textit{GL}_2}(P^{i})$. Recall that the $\textit{GL}_2$-equivariant Chow ring of $\mathbb{P}^i$ is generated by the Chern classes $\lambda_1, \lambda_2$ of the Hodge bundle and the first Chern class of $\mathcal{O}_{P^i}(-1)$, which we will call $t_i$, and the only relation is a polynomial $f_i(t_i,\lambda_1, \lambda_2)$ of degree $i+1$ (\cite[3.2, prop.6]{EdidinGraham} and the formula for projective bundles).

\begin{prop}\label{classi}
The class of $\Delta_{1,2i}$ in $\op{CH}^{1}_{\textit{GL}_2}(P^{2i})$ is always divisible by $2$. It is divisible by $p$ if and only if $2i-1$ is divisible by $p$.
\end{prop}
\begin{proof}
Consider the following commutative diagram:
\begin{center}
$\xymatrixcolsep{5pc}
\xymatrix{ (P^1)^{(2i-2)} \times P^1 \ar@{->}[d]^{\rho_1} \ar@{->}[r]^{i} & (P^1)^{2i} \ar@{->}[d]^{\rho_2}  \\ 
P^{2i-2}\times P^{1} \ar@{->}[r]^{\pi_{1,2i}} & P^{2i} } 
$

\end{center}

First note that the equivariant chow rings of all $\textit{GL}_2$-schemes involved are torsion-free, so we can make computations with rational numbers. The map $\rho_1$ is defined by $(f_1,\ldots,f_{2i-2},g) \rightarrow (f_1 \ldots f_{2i-2} , g)$, the map $\rho_2$ is defined by $(f_1,\ldots,f_{2i}) \rightarrow (f_1 \ldots f_{2i})$, and the map $i$ is defined by $(f_1,\ldots,f_{2i-2},g) \rightarrow (f_1, \ldots , f_{2i-2}, g, g)$. All the maps in the diagram are $\textit{GL}_{2}$-equivariant, $i$ is a closed immersion, $\rho_1, \rho_2$ are finite of degree respectively $(2i-2)!$ and $2i!$.

Now, the class of $\Delta_{1,2i}$ is the image of $1$ through $(\pi_{1,2i})_{*}$. Following the left side of the diagram we obtain $\left[ \Delta_{1,2i} \right] = \frac{1}{(2i-2)!} (\pi_{1,2i}\circ \rho_{1})_{*} (1)$. Consider now the right side of the diagram. The equivariant chow ring of $(P^{1})^{2i}$ is generated by all the different pullbacks of $t_1$, which we will call $y_1,\ldots , y_{2i}$, plus $\lambda_1$ and $\lambda_2$. It is easy to check that the pullback of $t_{2i}$ is $y_1 + \ldots + y_{2i}$, which by the projection formula, and by symmetry, implies that $(\rho_2)_{*}(y_j)=(2i-1)! t_{2i}$.

Using lemma \ref{Diag}, we see that $i_* (1) = y_i + y_{2i-1} + \lambda_1$. Its image is $2 (2i-1)! t_{2i} + 2i! \lambda_1$. By comparing the two formulae we obtain $$\left[ \Delta_{1,2i} \right]=\frac{1}{(2i-2)!}(2 (2i-1)! t_{2i} + 2i!\lambda_1)=2(2i-1)t_{2i}+(2i-1)2i \lambda_1$$ which implies the statement of the proposition. 
\end{proof}

The proposition above, together with the algebraic closure of our base field, will allow us to understand the pushforward $A^0_{\textit{GL}_2}(\Delta_{1,n})\rightarrow A^1_{\textit{GL}_2}(P^n)$. With the following proposition we establish a comparison between $A^0_{\textit{GL}_2}(\Delta_{1,n})$ and $A^0_{\textit{GL}_2}(\Delta_{1,n}\! \smallsetminus \! \Delta_{2,n})$.

As above, using the projective bundle formula and the computation of $A_{\textit{GL}_2}(\op{Spec}(k_0))$ \ref{proj},\ref{GLSL} we see that the equivariant Chow groups with coefficients $A^{\bull}_{\textit{GL}_2}(P^{2i})$ are just the tensor product of the usual equivariant Chow groups and $\op{H}^{\bull}(k_0)$, that is, they are generated as an $\op{H}^{\bull}(k_0)$-algebra by $\lambda_1, \lambda_2$ and $t_{2i}$, modulo a polynomial of degree $2i+1$ in $\lambda_1, \lambda_2, t_{2i}$ whose coefficients are in $\op{H}^{0}(k_0)=\mathbb{Z}/p\mathbb{Z}$.

Before moving on, we should do the following remark: suppose that $X$ is a scheme over $k_0$ such that there is an open subset $U$ of $X$ that is a smooth-Nisnevich covering of $\op{Spec}(k_0)$, that is, $U$ is smooth and has a rational point (see definition 3.2 of \cite{Pirisi}). Then by \cite[3.8]{Pirisi} the pullback $\op{H}^{\bull}(k_0)=\op{Inv}^{\bull}(\op{Spec}(k_0))\rightarrow \op{Inv}^{\bull}(U) = A^{0}(U)$ is injective.

 The pullback above factors through $A^{0}(X) \rightarrow A^{0}(U)$, so we can conclude that $\op{H}^{\bull}(k_0)$ maps injectively to $A^{0}(X)$. All of the equivariant approximations we will be using will have a smooth open subset with a rational point, so it makes sense to say that a group $A^{0}_{\textit{GL}_2}(T)$ is \emph{trivial} if it is equal to $\op{H}^{\bull}(k_0)$.

Next proposition describes the zero codimensional Chow group with coefficients $A^{0}_{\textit{GL}_2}(\Delta_{r,2i})$. 

\begin{prop}\label{oddeven}
If $r+1$ is divisible by $p$, the inclusion map $\Delta_{r,2i} \! \smallsetminus \! \Delta_{r+1,2i} \xrightarrow{j} \Delta_{r,2i}$ induces an isomorphism on $\op{A}^{0}_{\textit{GL}_2}$. If $r+1$ is not, the group $\op{A}^{0}_{\textit{GL}_2}(\Delta_{r,2i})$ is trivial, i.e. it is generated by $1$ as a (free) $\op{H}^{\bull}(k_0)$-module.
\end{prop}
\begin{proof}
 
  As $\op{A}^{0}_{\textit{GL}_2}(\Delta_{r,2i})$ is isomorphic to $\op{A}^{0}_{\textit{GL}_2}(\Delta_{r,2i} \! \smallsetminus \! \Delta_{r+2,2i})$ (because $\Delta_{r+2,2i}$ has codimension two in $\Delta_{r,2i}$) we can compute it using the following exact sequence:

 \begin{center}
 $ 0 \rightarrow \op{A}^{0}_{\textit{GL}_2}(\Delta_{r,2i} \! \smallsetminus \! \Delta_{r+2,2i}) \rightarrow \op{A}^{0}_{\textit{GL}_2}(\Delta_{r,2i} \! \smallsetminus \! \Delta_{r+1,2i}) \xrightarrow{\partial} \op{A}^{0}_{\textit{GL}_2}(\Delta_{r+1,2i} \! \smallsetminus \! \Delta_{r+2,2i})$
 \end{center}
 
 When $r+1$ is coprime to $p$, we want to prove that the kernel of $\partial$ is equal to $\op{H}^{\bull}(k_0)$. This will then imply that the image of $\op{A}^{0}_{\textit{GL}_2}(\Delta_{r,2i} \! \smallsetminus \! \Delta_{r+2,2i})$ must be equal to $\op{H}^{\bull}(k_0)$, and thus it must be trivial. When $r+1$ is divisible by $p$, we want to prove that $\partial$ is zero, so that the second arrow will be an isomorphism.
 
  The map $(P^{2i-2r} \! \smallsetminus \! \Delta_{2,2r}) \times P^{r} \xrightarrow{\pi} \Delta_{r,2i} \! \smallsetminus \! \Delta_{r+2,2i}$ yields the following commutative diagram with exact columns:

\begin{center}
$\xymatrixcolsep{3pc}
\xymatrix{ \op{A}^{0}_{\textit{GL}_2}((P^{2i-2r} {\! \smallsetminus \!} \Delta_{2,2i-2r}) \times P^r ) \ar@{->}[r]^{\pi_{*}} \ar@{->}[d] & \op{A}^{0}_{\textit{GL}_2}(\Delta_{r,2i} \! \smallsetminus \! \Delta_{r+2,2i} )   \ar@{->>}[d] \\
\op{A}^{0}_{\textit{GL}_2}((P^{2i-2r} \! \smallsetminus \! \Delta_{1,2i-2r}) \times P^r) \ar@{^{(}->>}[r]^{\pi_{*}} \ar@{->}[d]^{\partial_1} &  \op{A}^{0}_{\textit{GL}_2}(\Delta_{r,2i}\! \smallsetminus \! \Delta_{r+1,2i}) \ar@{->}[d]^{\partial}\\
 \op{A}^{0}_{\textit{GL}_2}((\Delta_{1,2i-2r} \! \smallsetminus \! \Delta_{2,2i-2r})\times P^r) \ar@{->}[r]^{\pi_{*}} & \op{A}^{0}_{\textit{GL}_{2}}(\Delta_{r+1,2i}\! \smallsetminus \! \Delta_{r+2,2i})  } $
\end{center}

 The second horizontal map is an isomorphism by \ref{univ},\ref{univiso} because $\pi_*$ is a universal homeomorphism when restricted to $\Delta_{r,2i} \! \smallsetminus \! \Delta_{r+1,2i}$.
 
  The kernel of $\partial_{1}$ is trivial because $\op{A}^{0}_{\textit{GL}_2}((P^{2i-2r} \! \smallsetminus \! \Delta_{2,2i-2r}) \times P^r)$ is trivial, as $A^{0}_{\textit{GL}_2}(P^{2i-2r} \times P^r)$ is trivial by the projective bundle formula \ref{proj} and $\Delta_{2,2i-2r} \times P^r$ has codimension $2$. 
 
  We claim that when $r+1$ is coprime to $p$ the third horizontal map is an isomorphism, implying that the kernel of $\partial$ must be trivial too, and when $r+1$ is divisible by $p$ the third horizontal map is zero, so that $\partial$ must be zero too.
 
 Let $\psi$ be the map from $(P^{2i-2r-2} \! \smallsetminus \! \Delta_{1,2i-2r-2})\times P^{r} \times P^1$ to $(P^{2i-2r-2} \! \smallsetminus \! \Delta_{1,2i-2r-2})\times P^{r+1}$ sending $(f,g,h)$ to $(f,gh)$. We have a commutative diagram:
 
\begin{center}
$\xymatrixcolsep{5pc}
\xymatrix{ (P^{2i-2r-2} \! \smallsetminus \! \Delta_{1,2i-2r-2}) \times P^1 \times P^{r}  \ar@{->}[r]^{\pi_1} \ar@{->}[d]^{\psi} & (\Delta_{1,2i -2r} \! \smallsetminus \! \Delta_{2,2i-2r}) \times P^r \ar@{->}[d]^{\pi}\\ 
(P^{2i-2r-2} \! \smallsetminus \! \Delta_{1,2i-2r-2})\times P^{r+1}  \ar@{->}[r]^{\pi_2} & \Delta_{r +1,2i} \! \smallsetminus \! \Delta_{r +2,2i}}$
\end{center}

Where $\pi_1$ and $\pi_2$ are defined respectively by $(f,g,h) \rightarrow (fg^2,h)$ and $(f,g) \rightarrow (fg^2)$. The maps $\pi_1$ and $\pi_2$ are universal homeomorphisms, so the pushforward maps $(\pi_1)_*,(\pi_2)_*$ are isomorphisms by \ref{univiso}. Then if we prove that $\psi_*$ is an isomorphism then $\pi_*$ will be an isomorphism too, and if $\psi_*$ is zero then $\pi_*$ will be zero too. Consider this last diagram:
 
\begin{center}
$\xymatrixcolsep{5pc}
\xymatrix{ (P^{2i-2r-2} \! \smallsetminus \! \Delta_{1,2i-2r-2})\times P^{r} \times P^1 \ar@{->}[dr]^{p_1} \ar@{->}[d]^{\psi}\\ 
(P^{2i-2r-2} \! \smallsetminus \! \Delta_{1,2i-2r-2})\times P^{r+1}  \ar@{->}[r]^{p_2} & P^{2i-2r-2} \! \smallsetminus \! \Delta_{1,2i-2r-2}}$
\end{center}
The pullbacks along $p_1$ and $p_2$ are both isomorphism by the projective bundle formula \ref{proj}, implying that the pullback of $\psi$ is an isomorphism. We have $\psi_* ( \psi^* \alpha) = \op{deg}(\psi) \alpha$ by the projection formula. Then as the degree of $\psi$ is $r + 1$, $\psi_*$ is zero if $p \mid r+1$ and an isomorphism zero otherwise.
 
\end{proof}

\begin{cor}\label{piinondue}
Let $p \neq 2$. If the class of $\Delta_{1,2i}$ is divisible by $p$ in $\op{CH}^{1}_{\textit{GL}_2}(P^{2i})$ then $\op{A}^{0}_{\textit{GL}_2}(P^{2i} \! \smallsetminus \! \Delta_{1,2i})$ is generated by $\langle 1, \alpha \rangle$, where $\alpha \neq 0$ is the invariant in degree $1$ corresponding to an equation for $\Delta_{1,2i}$. Otherwise $\op{A}^{0}_{\textit{GL}_2}(P^{2i} \! \smallsetminus \! \Delta_{1,2i})$ is trivial.
\end{cor}
\begin{proof}
The previous proposition shows that $\op{A}^{0}_{\textit{GL}_2}(\Delta_{1,2i})$ is always trivial for $p \neq 2$. Then applying the localization exact sequence to the inclusion $\Delta_{1,2i}\rightarrow P^{2i}$ yields the result.
\end{proof}

\begin{rem}\label{nonalg1}
The results \ref{oddeven}, \ref{piinondue} above do not require that $k_0$ be algebraically closed. However, this hypothesis will be fundamental in the next few steps.
\end{rem}

 From the next corollary on we will rely heavily on the algebraic closure of $k_0$, which is necessary to prove that the image of $i_* : A^{0}_{\textit{GL}_2}(\Delta_{1,2i})\rightarrow A^{0}_{\textit{GL}_2}(P^{2i})$ is zero.  In the next sections we will explore some ideas to get around this obstacle. Recall that for Chow rings with coefficients we will always refer to the grading coming from the cycle module as degree, and we will use the word codimension for the other.

\begin{cor}\label{piidue}
If $p=2$, then $\op{A}^{0}_{\textit{GL}_2}(P^{2i} \! \smallsetminus \! \Delta_{1,2i})$ is generated as a $\mathbb{F}_{2}$-module by $ \lbrace 1, \alpha_1, \ldots , \alpha_i \rbrace$, where the degree of $\alpha_i$ is $i$, and all the $\alpha_i$ are nonzero.
\end{cor}
\begin{proof}
As we are assuming that our base field is algebraically closed, we know that the Chow ring with coefficients $A^{\bull}_{\textit{GL}_2}(P^{2i})$ is contained in (cohomological) degree zero (because $\op{H}^{\bull}(k_0)=\op{H}^{0}(k_0)=\mathbb{Z}/p\mathbb{Z}$). Moreover the class of $\Delta_{1,2i}$ in $A^1_{\textit{GL}_2}(P^{2i})$ is divisible by $2$ due to proposition \ref{classi}, so the pushforward $A^0_{\textit{GL}_2}(\Delta_{1,2i}) \rightarrow A^1_{\textit{GL}_2}(P^{2i})$ must always be zero. 

Proposition \ref{oddeven} tells us that $A^0_{\textit{GL}_2}(\Delta_{1,2i})$ is isomorphic to $A^0_{\textit{GL}_2}(\Delta_{1,2i} \! \smallsetminus \! \Delta_{2,2i})$ which in turn is isomorphic to $A^{0}(P^{2i-2}\! \smallsetminus \! \Delta_{1,2i-2})$. Then we can setup an inductive reasoning using the localization exact sequence:

$$0 \rightarrow A^{0}_{\textit{GL}_2}(P^{2i}) \rightarrow A^{0}_{\textit{GL}_2}(P^{2i}\! \smallsetminus \! \Delta_{1,2i})  \rightarrow A^{0}_{\textit{GL}_2}(\Delta_{1,2i}) \rightarrow 0 $$

The first step is given by considering $\Delta_{1,2}$, which is universally homeomorphic to $P^1$, mapping to $P^2$. All the modules appearing in the exact sequences above are free, so the sequences all split. We can easily conclude.
\end{proof}

\begin{proof}[ Proof of Theorem \ref{evenmain}]

If $p = 2$ we start from corollary \ref{piidue}. The $\mathbb{G}_m$-bundle $$\left[\mathbb{A}^{4g+3}\! \smallsetminus \! \Delta/\textit{GL}_2\right] \rightarrow \left[P^{2g+2} \! \smallsetminus \! \Delta_{1,2g+2} /\textit{GL}_2 \right] $$ can be extended to a line bundle $$\mathscr{L}\xrightarrow{f} \left[P^{2g+2} \! \smallsetminus \! \Delta_{1,2g+2}\right].$$

 By taking a retraction $r$ we identify the equivariant chow groups with coefficients of $\mathscr{L}$ with those of $\left[ P^{2g+2} \! \smallsetminus \! \Delta_{1,2g+2} /\textit{GL}_2 \right]$. Then we can consider the exact sequence

\begin{center}
$0 \rightarrow \op{A}^{0}_{\textit{GL}_2}(P^{2g+2} \! \smallsetminus \! \Delta_{1,2g+2}) \xrightarrow{j^*\circ f^*} \op{A}^{0}_{\textit{GL}_2}(\mathbb{A}^{4g+3}\! \smallsetminus \! \Delta) \xrightarrow{\partial} \op{A}^{0}_{\textit{GL}_2}(P^{2g+2} \! \smallsetminus \! \Delta_{1,2g+2}) \xrightarrow{r \circ s_*} \op{A}^{1}_{\textit{GL}_2}(P^{2g+2} \! \smallsetminus \! \Delta_{1,2g+2})$
\end{center}

Where $s$ is the zero section of $\mathscr{L}$ and $j$ is the open immersion $\left[\mathbb{A}^{4g+3}\! \smallsetminus \! \Delta/\textit{GL}_2\right] \xrightarrow{j} \mathscr{L}$. We want to understand the kernel of the map $r\circ s_*$. Using definition \ref{C1def} we can identify $r\circ s_*$ with the first Chern class of the line bundle $\mathscr{L}$, which by \cite[3.2]{FulghesuEdidin} is equal to $g\lambda_1 - t_{2g+2}$. As $g$ is even we have $g \lambda=0$ and our claim boils down to understanding whether the products $t_{2g+2} \alpha_j$ are zero or not.

Recall that by \ref{piidue} we have $A^0_{\textit{GL}_2}(P^{2i} \! \smallsetminus \! \Delta_{1,2i})= \langle 1, \alpha_1, \ldots , \alpha_i \rangle$, so by writing $\alpha_j$ we mean the generator in degree $j$ for some given $i$, which should be clear from context.

We will proceed by induction on $i$. First we take a look at the products in $\op{A}^{\bull}_{\textit{GL}_2}(P^2 \! \smallsetminus \! \Delta_{1,2})$. The second part of the localization exact sequence for $\Delta_{1,2} \rightarrow P^2$ reads:

\begin{center}
$0 \rightarrow \op{A}^{1}_{\textit{GL}_2}(P^2) \rightarrow \op{A}^{1}_{\textit{GL}_2}(P^2 \! \smallsetminus \! \Delta_{1,2}) \xrightarrow{\partial} \op{A}^{1}_{\textit{GL}_2}(\Delta_{1,2})$
\end{center}

We need to understand what $t_2 \alpha_1$ is. By the compatibility of Chern classes and boundary maps \ref{C1}, (v) we know that $\partial(t_2 \alpha_1) = \partial( c_1(\mathcal{O}_{P^2}(-1))(\alpha_1))= c_1(i^*\mathcal{O}_{P^2}(-1))(\partial(\alpha_1))=c_1(i^*\mathcal{O}_{P^2}(-1))(1)$. As the pullback of $\mathcal{O}_{P^2}(-1)$ through $P^1 \xrightarrow{\pi_{1,2}} \Delta_{1,2} \rightarrow P^2$ is equal to $\mathcal{O}_{P^1}(-1)^2$, we see that $\partial(t_2 \alpha_1)=0$. Then by the exact sequence above $t_2 \alpha_1$ must be the image of some $\gamma \in \op{A}^{1}_{\textit{GL}_2}(P^2)$, but there are no element of positive degree in $A^{\bull}_{\textit{GL}_2}(P^2)$ when $k_0$ is algebraically closed, so $t_2 \alpha_1 = 0$.

 Suppose by induction that for $2i < 2g+2$ we know that $t_{2i} \alpha_j \in \op{A}^{0}_{\textit{GL}_2}(P^{2i} \! \smallsetminus \! \Delta_{1,2i})$ is equal to zero if and only if $j = i$. We already know that for $i =1$, giving us the base for the induction. Again, we consider the exact sequence

\begin{center}
$0 \rightarrow \op{A}^{1}_{\textit{GL}_2}(P^{2g + 2} \! \smallsetminus \! \Delta_{2,2g+2}) \rightarrow \op{A}^{1}_{\textit{GL}_2}(P^{2g+2} \! \smallsetminus \! \Delta_{1,2g+2}) \xrightarrow{\partial} \op{A}^{1}_{\textit{GL}_2}(\Delta_{1,2g+2} \! \smallsetminus \! \Delta_{2,2g+2})$
\end{center}
we can see, for example by looking at the proof of \ref{classi}, that the pullback of $\mathcal{O}_{P^{2g+2}}(-1)$ to $\Delta_{1,2g+2}$ is equal to $\mathcal{O}_{P^{2g}}(-1) \otimes \mathcal{O}_{P^1}(-1)^2$, whose Chern class is the same as $c_1(\mathcal{O}_{P^{2g}}(-1))$ as we are working modulo $2$.

 Then by \ref{C1}, (v) we have of $\partial (t_{2g+2} \alpha_j) = (\pi_{1,2g+2})_* t_{2g} \alpha_{j-1}$, where we consider $\alpha_0 = 1$, showing that $t_{2g+2} \alpha_j\neq 0$. This immediately implies the thesis for $j < g+1$.

 Now consider the case where $j=i=g+1$. As $$\partial (t_{2g+2} \alpha_{g+1}) = (\pi_{1,2g+2})_* t_{2g} \alpha_{g} = 0$$ the element $t_{2g+2} \alpha_{g+1}$ must be the image of an element of $\op{A}^{1}_{\textit{GL}_2}(P^{2g + 2} \! \smallsetminus \! \Delta_{2,2g+2})$. The elements of $\op{A}^{1}_{\textit{GL}_2}(P^{2g + 2} \! \smallsetminus \! \Delta_{2,2g+2})$ can have degrees only up to one plus the maximum degree of an element of $\op{A}^{0}_{\textit{GL}_2}(\Delta_{2,2g+2})$, again due to the localization exact sequence and the fact that the Chow ring with coefficients $A^{\bull}_{\textit{GL}_2}(P^{2g+2})$ is contained in degree zero. Then by proposition \ref{oddeven} their degree is equal or less than one, and $t_{2g+2} \alpha_{g+1}$ must be zero, concluding the proof.

If $p \neq 2$, starting from corollary \ref{piinondue} we only have to possibly check that $g\lambda_1 - t_{2g+2}\alpha_1$ is not $0$, which can be done exactly as above. The explicit result for $p \neq 2$ can be obtained by looking at whether the class of $\Delta_{1,2g + 2}$ is divisible by $p$ in the equivariant Picard group of $P^{2g+2}$, which can be easily done using proposition \ref{classi}.

\end{proof}

\begin{rem}
For $g=2$, we understand the multiplicative structure of $\op{Inv}^{\bull}(\mathscr{M}_2)$ almost completely. We have $\alpha_4^2 = \alpha_4 \alpha_3 = \alpha_4 \alpha_2 = \alpha_4 \alpha_1= 0$ as there are no elements of degree higher than $\alpha_4$, and similarly $\alpha_3^2 = \alpha_3 \alpha_2 =  \alpha_3 \alpha_1= 0$ as these elements are pullbacks from $\op{Inv}^{\bull}(\left[ P^6 \! \smallsetminus \! \Delta_{1,6} / \textit{GL}_2 \right])$ and we can apply the same reasoning. 

The squares $\alpha_1^2, \alpha_2^2$ are both zero, as the second is of degree $4$ and there are no elements of degree $4$ in $\op{Inv}^{\bull}(\left[ P^6 \! \smallsetminus \! \Delta_{1,6} / \textit{GL}_2 \right])$, and the first is represented by the square of an element $\tilde{\alpha_1} \in \op{H}^1(k(P^6))=k(P^6)^*/(k(P^6)^*)^2$ which is equal to the element $\lbrace -1 \rbrace \tilde{\alpha_1} \in \op{H}^2(k(P^6))$ which is zero as $k$ contains a square root of $-1$. The product $\alpha_1 \alpha_2$ may be either equal to zero or to $\alpha_3$.

In general we have no instruments to understand the multiplicative structure of $\op{Inv}^{\bull}(\mathscr{H}_g)$, and considerations like the above are less and less useful as the number of generators and possible degrees grows. The problem seems to be that it is very difficult to keep track of how the products of our elements behave when using the localization exact sequence, and in fact in most computations (that the author knows of) of classical cohomological invariants the multiplicative structure stems from an explicit \emph{a priori} description of the invariants.
\end{rem}
\section{The non algebraically closed case}

In this section we obtain a partial result on the cohomological invariants of $\mathscr{M}_2$ for a general base field. This should give an idea of the inherent problems that arise when we have nontrivial elements of positive degree in the cohomological invariants of our base field.

 Note that this will happen even for an algebraically closed field if we are considering quotients by groups that are not special, making the development of techniques and ideas to treat these type of problems crucial for the future development of the theory.

\begin{thm}\label{M2general}
Suppose that the characteristic of $k_0$ is different from $2,3$.

 Let $p$ be a prime different from $2$. Then the cohomological invariants of $\mathscr{H}_g$ are nontrivial if and only if $4g + 1$ is divisible by $p$. In this case they are freely generated as a $\op{H}^{\bull}(\op{Spec}(k))$-module by $1$ and a single nonzero invariant in degree one.

Let $p$ be equal to $2$. Then the cohomological invariants of $\mathscr{M}_2$ fit into the following exact sequence of $\op{H}^{\bull}(\op{Spec}(k))$-modules 
$$ 0 \rightarrow M \rightarrow \op{Inv}^{\bull}(\mathscr{M}_2) \rightarrow K \rightarrow 0$$
 where $M$ is freely generated by $1$ and elements $\alpha_1,\alpha_2,\alpha_3$ of respective degrees $1,2,3$ and $K$ is isomorphic to a submodule of  $\op{H}^{\bull}(\op{Spec}(k))$, shifted in degree by $4$.
\end{thm}

The first statement is a direct consequence of remark \ref{nonalg1}. The part of the proof of \ref{evenmain} where $p\neq 2$ can be carried out exactly in the same way. The case $p=2$ will require some work, and in the rest of the section we always work in this case.

There are two points that we need to prove for the machinery we used in the previous section to work:
\begin{enumerate}
\item The pullback through the map $\Delta_{r,2i} \! \smallsetminus \! \Delta_{r+1,2i} \rightarrow \Delta_{r,2i}$ must induce an isomorphism on $A^0$ for $i \leq 3$, as in proposition \ref{oddeven}.

\item The pushforward through the map $\Delta_{1,2i} \rightarrow P^{2i}$ must be zero for $i \leq 3$, as in corollary \ref{piidue}.
\end{enumerate}

The first point is again implied by remark \ref{nonalg1}. The proof of the second point in the previous section completely depends on $k_0$ being algebraically closed, so we will have to think of something new. We begin by reducing to the non-equivariant case when $i=3$.

 The main idea of the section is that we can use the properties of Chern classes to construct an element $f$ in $A^{\bull}_{\textit{GL}_2}(P^i)$ (resp. $A^{\bull}(P^i)$) such that $f$ annihilates the image of $A^{0}_{\textit{GL}_2}(\Delta_{1,2i})$ (resp. $A^{0}(\Delta_{1,2i})$) but multiplication by $f$ is injective on $A_{\textit{GL}_2}^1(P^{2i})$ (resp. $A^1(P^{2i})$).

\begin{lem}\label{lemnonalg}
The map $A^0_{\textit{GL}_2}(\Delta_{1,6})\rightarrow A^1_{\textit{GL}_2}(P^6)$ is zero if and only if the map $A^0(\Delta_{1,6}) \rightarrow A^1(P^{6})$ is zero.
\end{lem}
\begin{proof}
One arrow is trivial: the equivariant groups for $P^6$ map surjectively to the non-equivariant groups and the assignment is functorial, so if the equivariant map is trivial the same must be true for the non-equivariant map.

We now remove $\Delta_{2,6}$ from both sides, so that we are reduced to considering the map $(P^4\! \smallsetminus \! \Delta_{1,4})\times P^1 \rightarrow P^6 \! \smallsetminus \! \Delta_{2,6}$. All elements in $A^0_{\textit{GL}_2}((P^4\! \smallsetminus \! \Delta_{1,4})\times P^1)$ are pullbacks through the first projection. Following \cite[p.638]{Vistoli} and using the fact that we are working modulo $2$, we see that an element $\alpha \in A^0_{\textit{GL}_2}(P^4 \! \smallsetminus \! \Delta_{1,4})$ satisfies the equation $(t_4^5 + t_4^3 \lambda_1^2)\alpha = 0$. The pullback $\alpha \in A^0((P^4\! \smallsetminus \! \Delta_{1,4})\times P^1)$ must then satisfy the same equation.

Recall now that modulo two the pullback of $\mathcal{O}_{P^6}(-1)$ is equal to $\mathcal{O}_{P^4}(-1)$. As $\lambda_1$ is also the pullback of the corresponding equivariant line bundle on $P^6$, by the projection formula we see that the image of $\alpha$ must satisfy the same equation. As $i_*(\alpha)$ is an element of $A^1_{\textit{GL}_2}(P^6)$ we can write $i_*(\alpha)=\lambda_1 \cdot a + t_6 \cdot b$ with $a,b \in H^{\bull}(\op{spec}(k))$. Then we have $(t_6^5 + t_6^3 \lambda_1^2)( \lambda_1 \cdot a + t_6 \cdot b)=0$ in $A^6_{\textit{GL}_2}(P^6 \! \smallsetminus \! \Delta_{2,6})$.

Suppose that we know the result in the non-equivariant case, that is, we know that $b = 0$. We want to show that $(t_6^5 + t_6^3 \lambda_1^2)\lambda_1 \cdot a$ belongs to the image of $A^4_{\textit{GL}_2}(\Delta_{2,6})$ if and only if $a=0$. Recall that $\Delta_{2,6}$ can be seen as the union of $(P^2 \! \smallsetminus \! \Delta_{1,2})\times P^2$ and $\Delta_{3,6}$. The elements in $A^{\bull}_{\textit{GL}_2}(\Delta_{2,6})$ are sums of elements of three types: those that come from $P^2 \times P^2$, those that come from $\Delta_{3,6}$ (which is universally homeomorphic to $P^3$) and the elements of $A^{\bull}_{\textit{GL}_2}((P^2 \! \smallsetminus \! \Delta_{1,2})\times P^2)$ that are ramified on $\Delta_{1,2} \times P_2$ but unramified on $\Delta_{3,6}$.

  Using again the computations in \cite{Vistoli} we see that the first two images form the ideal $(t_6^6 + t_6^5 \lambda_1 + t_6^4 \lambda^2_1 + t_6^3 \lambda_1^3)$. For the latter, the computation reduces to finding out the kernel of the pushforward $A^{\bull}_{\textit{GL}_2}(P^1 \times P^2) \rightarrow A^{\bull}(P^3)$. Using the fact that the map is finite of degree $3$ one sees that if we write $t = c_1(\mathcal{O}_{P^1}(-1)), s = c_1(\mathcal{O}_{P^2}(-1))$ the kernel is generated as a $A^{\bull}_{\textit{GL}_2}(\op{Spec}(k))$-module by $1,s,st$. Then any element in codimension $4$ belonging to the kernel of our pushforward can be written down as a sum $\lambda_1^2 a_1 + \lambda_2 a_2$, and the same must hold for any element in $A^4_{\textit{GL}_2}(\Delta_{2,6})$ belonging to the third type (up to elements coming from $P^2 \times P^2$). By the projection formula we can conclude that the image of $A^4_{\textit{GL}_2}(\Delta_{2,6})$ must be contained in the ideal $(t_6^6 + t_6^5 \lambda_1 + t_6^4 \lambda^2_1 + t_6^3 \lambda_1^3, \lambda_1^2, \lambda_2)$, which does not contain $(t_6^5 + t_6^3 \lambda_1^2)\lambda_1 \cdot a$ unless $a=0$.
\end{proof}

Of course the same trick will not work on the non-equivariant case, as the relation would be $t_6^6 \cdot a = 0$ and $\Delta_{2,6}$ contains rational points. We will have to dirty our hands and work at cycle level. Recall that the  the first Chern class of a line bundle $L$ can be defined on cycles by choosing a coordination for $L$.

\begin{lem}
Let $E \rightarrow X$ be a line bundle that is isomorphic to $L \otimes W^{\otimes p}$ for some line bundles $L,W$. Let $\tau_L, \tau_W$ be coordinations respectively for $L$ and $W$, and consider the coordination $\tau'=\tau_L \cup \tau_W$ for $E$. Additionally let $\tau''$ be the coordination for $L$ that is set theoretically $\tau_L \cup \tau_W$, with the trivializations induced by those of $\tau_L$. Then for all $\alpha \in C^0(X)$ we have $c_{1,\tau'}(E)(\alpha) = c_{1,\tau''}(L)(\alpha)$.
\end{lem}
\begin{proof}
Fix the additional data of an open covering $U=\sqcup U_i$ of $X$ such that $E, L, W$ are all trivial on $U$. Then the three bundles are respectively described by a choice of coordinate change elements $\alpha_{i,j} \in \mathcal{O}^*(U_i \times_X U_j)$ for each couple $(i,j)$, and we can take $\alpha_{i,j,E}=\alpha_{i,j,L} \cdot \alpha_{i,j,W}^p$. It can be seen directly as in \ref{C1}, (iv) that given a compatible choice of a trivialization and coordination for $E$ the Chern class $c_{1,\tau'}(E)(\alpha)$ can be decomposed (not uniquely) as the sum of $c_{1,\tau_L \cup \tau_W}(X \times \mathbb{A}^1)(\alpha)$ and a function that is linear in the coordinate change elements $\alpha_{i,j,E }\in \mathcal{O}^*(U_i \times_X U_j)$.

 As the elements $\alpha_{i,j,E}$ satisfy $\alpha_{i,j,E}=\alpha_{i,j,L}\cdot \alpha_{i,j,W}^p$, and a $p$-th power is zero in Galois cohomology with coefficients in $\mu_p$, we can conclude.
\end{proof}

Consider now the line bundle $\mathcal{O}(-1)$ on $P^{2i}$, with the coordination given by repeatedly taking the hyperplane at infinity. We think of $P^{2i}$ as the space of forms of degree $2i$ up to scalars, by $$(x_0 : x_1 : \ldots : x_{2i}) \rightarrow \left[ x_0 s_0^{2i} + \alpha_1 s_0^{2i-1}s_1 + \ldots + x_{2i}s_1^{2i} \right]$$
so we can think of the $l$-th element of the coordination as imposing that the first $l$ coefficients are zero. We will name this coordination $\sigma$.

The pullback of $\sigma$ to $P^{2i-2}\times P^1$ looks somewhat strange. We want to think of the map $P^{2i-2}\times P^1 \rightarrow P^{2i}$ as $(f,g) \rightarrow fg^2$, so we are imposing that first $l$ coefficients of $fg^2$ must be zero. One then easily checks that the pullback of the coordination above is in the form

$$ P^{2i-2}\times P^1 \supset H_0 \times P^1 \cup P^{2i-2} \times \lbrace p \rbrace \supset H_1 \times P^1 \cup P^{2i-2} \times \lbrace p\rbrace \supset H_2 \times P^1 \cup H_0 \times \lbrace p \rbrace \supset \ldots $$

Where $H_0 \supset H_1 \supset H_2 \ldots $ is the same coordination as above on $P^{2i-2}$, $p$ corresponds to the point $(0:1) \in P^1$ and the terms in the form $X \times \lbrace p \rbrace$ each repeat twice as imposing the condition on $g$ kills two coefficients at a time. We will name this second coordination $\tau$.

Note now that the pullback of $\mathcal{O}(-1)$ to $P^{2i-2}\times P^1$ is $\mathcal{O}_{P^{2i-2}}(-1) \otimes \mathcal{O}_{P^1}(-1)^2$. Then by the lemma above taking the Chern class of the pullback of $\mathcal{O}(-1)$ with respect to $\tau$ is the same as taking $c_{1,\tau}(\mathcal{O}_{P^{2i-2}}(-1))$.

Note that the considerations above are still perfectly true for $P^n$ when $n$ is odd. We just never had to consider this case. We will need to do it in the following lemma and proposition.

\begin{lem}
Let $n$ be a positive integer (possibly odd). Let $\alpha$ be an element in $C^0(P^{n-2})$, possibly ramified only at $\Delta_{1,n-2}$. Let $\alpha'$ be its pullback to $C^0(P^{n-2}\times P^1)$. Then $c_{1,\tau}(\mathcal{O}_{P^{n-2}}(-1))(\alpha')$ belongs to $C^0(H_0 \times P^1)$, and it is the pullback of an element $\beta \in C^0(H)$, possibly ramified only at $H \cap \Delta_{1,n-2}$.
\end{lem}
\begin{proof}
We follow the local computation in the proof in \ref{C1}, (iv). There are two main differences. The first is that the first step of the coordination we are working with is not an irreducible divisor, but rather the union of two irreducible divisors. This will not pose a problem as it's easy to explicitly check that the part of the map coming from $P^{2i-2}\times \lbrace p \rbrace$ does not contribute to $c_{1,\tau}(\mathcal{O}_{P^{n-2}}(-1))(\alpha')$. The second difference is that our element $\alpha'$ is not unramified. The hypothesis that $\alpha'$ should be not ramified is only needed at the end, and the reasoning in \ref{C1} (iv) only uses the fact that it is not ramified at the generic point of our irreducible divisor, which is true in the case of $H$. Following this we see that $c_{1,\tau}(\mathcal{O}_{P^{n-2}}(-1))(\alpha')$ is equal to $c_{1,\sigma'}(\mathcal{O}_{P^{n-2}}(-1))(\alpha')$, where $\sigma'$ denotes the coordination on $P^{n-2}\times P^1$ obtain by pulling back $\sigma$ through the first projection. Then the claim follows from the compatibility of the first Chern class with pullbacks and differentials.
\end{proof}

\begin{prop}
Let $\alpha$ be as above. We have $c_{1,\tau}(\mathcal{O}_{P^{n-2}}(-1))^{n-1}(\alpha')=0$.
\end{prop}
\begin{proof}
Consider the element $c_{1,\tau}(\mathcal{O}_{P^{n-2}}(-1))(\alpha')$. By the lemma above it is the pushforward of some element $\beta$ in $H \times P^1$. We have $H \simeq P^{n-3}$, and the isomorphism sends $H \cap \Delta_{1,n-2}$ to $\Delta_{1,n-3}$. By the projection formula, the Chern class $c_{1,\tau}(\mathcal{O}_{P^{n-2}}(-1))(\alpha')$ can be computed on $H \times P^1 \simeq P^{n-3} \times P^1$. It's easy to see that we are in the same situation as the previous lemma. By applying this reasoning $n-1$ times we will eventually be in the situation where $H$ is just a point and we get zero.
\end{proof}

\begin{cor}\label{cornonalg}
The map $A^{0}(\Delta_{1,2i})\xrightarrow{i} A^1(P^{2i})$ is zero for all $i$.
\end{cor}
\begin{proof}
Consider an element $\gamma \in A^{0}(\Delta_{1,2i})$. It can be seen as the pushforward of an element $\alpha' \in C^0(P^{2i-2}\times P^1)$, satisfying the hypothesis of the lemma above. By the previous proposition we know that $c_{1,\tau}(\mathcal{O}_{P^{2i-2}}(-1))^{2i-1}(\alpha')=0$. By the projection formula this tells us that $t_{2i}^{2i-1} i_* \alpha' = 0$, but the only element in $A^{1}(P^{2i})$ that satisfies this equation is $0$, so $i_* \alpha' = 0$.
\end{proof}

\begin{proof}[ Proof of Theorem \ref{M2general}]
Once we put together the results we accumulated, we will prove the theorem in the same way as we did for \ref{evenmain}.

We have mainly concerned ourselves with $P^6$ up to now, so we should begin by tackling the lower dimension cases. Then the results in this section will allow us to conclude easily.

\begin{enumerate}

\item The maps $$A^0_{\textit{GL}_2}(\Delta_{1,2})\xrightarrow{i_*} A^1_{\textit{GL}_2}(P^2), \quad A^0_{\textit{GL}_2}(\Delta_{1,4})\rightarrow A^1_{\textit{GL}_2}(P^4)$$ are both zero.

 The first statement is due to the projection formula. All elements in $A^{0}_{\textit{GL}_2}(\Delta_{1,2})=\op{H}^{\bull}(k_0)$ are pullbacks from $A^{0}_{\textit{GL}_2}(P^2)$ so we have $i_*(\alpha)=i_*(i^*\alpha)=i_*(1)\alpha$. As $i_*(1)=\left[ \Delta_{1,2} \right] =0$ by \ref{classi}, we can conclude.
 
  To check the second statement, note that by the previous point and point (iii) below we have that $A^0_{\textit{GL}_2}(\Delta_{1,4})=\langle 1, \alpha \rangle$ as a free $H^*(\op{Spec}(k))$-module, where $\alpha$ is an element of degree $1$, coming from the cohomological invariants of $\left[ P^2 \!\smallsetminus\!\Delta_{1,2} /GL_2 \right]$. The image of $1$ is zero as the class of $\Delta_{1,4}$ in $A^1_{\textit{GL}_2}(P^4)$ is even.

Consider now the boundary map $\partial: A^0_{\textit{GL}_2}((P^2\!\smallsetminus\!\Delta_{1,2})\times P^1) \rightarrow A^{0}_{\textit{GL}_2}(\Delta_{1,2}\times P^1)$; we have $\partial (t_2 \cdot \alpha) = 0$, implying that $t_2 \alpha$ comes from $A^{0}_{\textit{GL}_2}(P^2\times P^1)$. Moreover, by the compatibility with flat pullback \ref{C1}, (ii) we know that $t_2 \alpha$ is a pullback from $A^1_{\textit{GL}_2}(P^2)$, and thus it can be written as $a t_2 + b \lambda_1$ with $a,b \in \op{H}^{\bull}(k_0)$. Using again the projection formula and proposition \ref{classi} we see we that the image of an element in this form must be zero, so we must have $t_4 i_* \alpha=0$. By the structure of the Chow groups with coefficient of a projective bundle this implies $i_* \alpha = 0$.

\item The map $A^0_{\textit{GL}_2}(\Delta_{1,6})\rightarrow A^1_{\textit{GL}_2}(P^6)$ is zero. This is due to \ref{lemnonalg} and \ref{cornonalg}.

\item The pullback $A^0_{\textit{GL}_2}(\Delta_{1,2i}) \rightarrow A^0_{\textit{GL}_2}(\Delta_{1,2i} \! \smallsetminus \! \Delta_{2,2i})$ is an isomorphism. This is remark \ref{nonalg1}.

\item Using the localization exact sequence, the points above easily imply that $A^0_{\textit{GL}_2}(P^6 \! \smallsetminus \! \Delta_{1,6})$ is freely generated as an $H^*(\op{Spec}(k))$-module by $1$ and elements $\alpha_1,\alpha_2,\alpha_3$ of degree respectively $1,2,3$. This can be done exactly as in \ref{piidue}.

 As in the previous section all that is left is to understand the kernel of $$c_1(\mathcal{O}_{P^6}(-1)):A^0_{\textit{GL}_2}(P^6 \! \smallsetminus \! \Delta_{1,6}) \rightarrow A^0_{\textit{GL}_2}(P^6 \! \smallsetminus \! \Delta_{1,6}).$$ we can proceed as in the previous section to prove by induction that the map is injective on the submodule generated by $1,\alpha_1,\alpha_2$. Unfortunately the reasoning we used before to prove that $\alpha_3$ must belong to the kernel of $c_1(\mathcal{O}_{P^6}(-1))$ no longer works, as it relied heavily on the algebraic closure of $k$, so we have to add the unspecified module $K\simeq \op{Ker}(c_1)\cap\langle \alpha_3 \rangle$ to our final result.
\end{enumerate}
\end{proof}


\begin{thebibliography}{{Sta}15}

\bibitem[AV04]{ArsieVistoli}
Alessandro Arsie and Angelo Vistoli, \emph{Stacks of cyclic covers of
  projective spaces}, Compositio Mathematica 140, 647-666 (2004).

\bibitem[EG96]{EdidinGraham}
Dan Edidin and William Graham, \emph{Equivariant intersection theory}, Invent.
  Math. 131 (1996).

\bibitem[EKM08]{ElmanKarpenkoMerkurjev}
Richard Elman, Nikita Karpenko, and Alexander Merkurjev, \emph{The algebraic
  and geometric theory of quadratic forms}, American Mathematical Society,
  2008.

\bibitem[FE09]{FulghesuEdidin}
Damiano Fulghesu and Dan Edidin, \emph{The integral chow ring of the stack of
  hyperelliptic curves of even genus}, Math Research Letters {v.16}
  (2009), no.~1, 27-40.

\bibitem[Gui08]{Guillot}
Pierre Guillot, \emph{Geometric methods for cohomological invariants},
  Documenta Mathematica {vol.12 521-545.} (2008).

\bibitem[Pir15a]{Pirisi15}
Roberto Pirisi, \emph{Cohomological invariants of algebraic curves}, Ph.D.
  thesis, Scuola Normale Superiore di Pisa, \url{http://www.math.ubc.ca/~rpirisi/thesis.pdf}, 2015.

\bibitem[Pir15b]{Pirisi}
Roberto Pirisi, \emph{Cohomological invariants of algebraic stacks}, ar{X}iv
  math.AG:1412.0554v2 (2015).

\bibitem[Ros96]{Rost}
Markus Rost, \emph{Chow groups with coefficients}, Documenta Mathematica
  {vol.1,319-393} (1996).


\bibitem[{Sta}15]{StacksProject}
The {Stacks Project Authors}, \emph{Stacks project},
  \url{http://stacks.math.columbia.edu}, 2015.

\bibitem[Tot99]{Totaro}
Burt Totaro, \emph{The chow ring of a classifying space}, Proceedings of
  Symposia in Pure Mathematics {v. 67} (1999), 249-281.

\bibitem[Vez00]{Vezzosi}
Gabriele Vezzosi, \emph{On the chow ring of the classifying stack of $\mbox{PGL}(3)$},
  Journal f\"ur die reine und angewandte Mathematik (Crelle) {No.
  523,1-54} (2000).

\bibitem[Vis96]{Vistoli}
Angelo Vistoli, \emph{The chow ring of $\mathscr{M}_2$, appendix to equivariant
  intersection theory}, Invent. Math {131} (1996).

\bibitem[VM06]{MolinaVistoli}
Angelo Vistoli and Alberto~Luis Molina, \emph{On the chow ring of classifying
  spaces for classical groups}, Rend. Sem. Mat. Univ. Padova {v.116}
  (2006).

\end{thebibliography}
\end{document}